# Students' experiences of learning Sciences during the Covid-19 pandemic and their suggestions for the next day: A Greek University Department case study

**Ioannis Rizos** & **Nikolaos Gkrekas**

Department of Mathematics, University of Thessaly, Lamia 35100, Greece

**Abstract**

Around the globe Covid-19 pandemic has influenced not only the education, but also our everyday life, among other aspects. In Greece, distance learning started to get in use widely in tertiary education, since the first national lockdown was announced and reshaped education in many ways. In the University of Thessaly, in the Department of Mathematics, undergraduate students opt a lot of different courses to attend and due to the Covid-19 crisis all of them are taught via web platforms. Some of the most significant theoretical subjects are «Calculus» (with applications in Science and Mechanics), «Physics» (Classical Mechanics) and «Philosophy of Science». In addition, some other applied subjects are «Programming Languages» and «Digital Technologies in Mathematics Education» and the impact the above have, generally in education and society itself. In this paper, we describe the different ways students have reacted regarding learning Sciences in a distance learning environment. We split our case study in two parts. The first one is about the way students experience e-learning and the second one is about their suggestions for the next day. Integrating e-questionnaires and interviews and taking into account parameters like economic factors and the permanent residence issue, we asked the students about their preferences among face-to-face learning, distance learning and a blended model. Remarks about the academic life and the possible ways of taking the extra step, after the Covid-19 crisis ceases to exist, are made.

# 1. Introduction

As the first months of 2020 passed until we accepted as a society that the coronavirus pandemic was real, it became apparent that in addition to doctors, nurses and many other brave people self-sacrifice to support our daily survival, academics had the opportunity and the moral duty to start conducting in-depth studies of current events (Jandrić et al. 2020). Thus, as many colleagues did, in mid-March 2020 we wrote a short article (Rizos 2020) about the role of mathematical models in epidemiology (the law of exponential growth, the Malthusian growth model, probabilistic models etc.) in accordance with the philosophical background of these models and we made some pedagogical considerations.

Then, on the occasion of that article, we organized an e-seminar in the Department of Mathematics in the University of Thessaly where we discussed with our students the extraordinary educational situation due to the coronavirus pandemic, the importance of Mathematics and the need to use digital technologies in distance Science and Technology teaching. Since then we have been in constant contact with them trying to face the challenges together, to keep alive their interest in learning and to improve the educational process. In the spring semester 2020-2021, after a whole year of distance learning at the University, we asked our students to share with us their personal experiences and opinions on distance education in Science and Technology and to submit their suggestions for the next day. The context and results of this research are presented in the following paragraphs.

# 2. The impact of the pandemic in education

During this period, education seems to be in the middle of a crisis. Roughly 258 million children, adolescents and youth were out of school in 2018 (UNESCO 2019). When Covid-19 pandemic started, this number became even larger due to lack of equipment and resources especially in the least developed countries (Mishra et al. 2020). More than 1.5 billion children and youth from 188 countries have been badly influenced by the sudden "lock" in schools and institutions (United Nations 2020).

At the same time, the pandemic really appears to have changed the established face-to-face education model. Distance education using ICT became from a very rare choice (usually in supplementary/ further education programs), to a common practice. Especially higher education changed radically. The way of teaching and learning Sciences and Technology has been reformed and the use of digital technologies has been upgraded. The widespread use of technology has profoundly affected both teachers and students in the field of Mathematics as well. On the one hand, some older professors cannot comprehend the new digital platforms and they seem to be unable to conduct a complete lecture based on them. Many math teachers prefer to teach by the old-fashioned and conventional way and the current situation looks like a dead end to them (Naidoo 2020). On the other hand, the students seem to be able to really gain a lot of knowledge and skills from this experience.

The pandemic has had a serious impact not only on the relationship between students and their teachers, but also between students themselves. The academic life has vastly changed and the way people interact with each other has been transformed. The interactivity of each session is one of the most important parts of the learning procedure. In these times, it is visible that human contact and interactivity have faded, because of the exclusive use of digital technologies.

As in all universities across the world, Greek universities have been through a variety of challenges during the pandemic. Challenges like the use of web platforms and the adaption to different time limits and schedules. In the following paragraphs we describe a research we conducted in order to carefully examine these challenges in undergraduate students' lives and we discuss their suggestions for the post-pandemic period.

## 3. The situation in a Greek University Department – An a priori analysis

Since the beginning of 2020, due to the Covid-19 pandemic, the social reality and economic activity of all countries in the world have been significantly affected. Thus, in Greece, where thousands of people have lost their lives and the economy has been hit hard, there is an obvious sense of *uncertainty*. Students could not remain unaffected by this situation. The institutions in which they study were closed, students moved back to their parents' homes and therefore lost their independence and a large part of their interactions with their peers.

Greek universities quickly transformed the teaching model they followed from face-to-face teaching in classroom, to distance learning using ICT in order to ensure the continuity in educational service. Platforms such as MS-Teams and WebEx as well as courses management systems like eClass and Moodle were used, covering both synchronous and asynchronous distance education, while online guidance and technical support on the tools were provided to students and instructors.

The Department of Mathematics of the University of Thessaly faces the same challenges as those of other university departments in Greece and in the rest European countries. Since March 11, 2020, when the first national lockdown was imposed, all courses have been taught online. However, there are some peculiarities. Although the University of Thessaly is one of the biggest and most respectful institutions in Greece and according to [Ranking Web of Universities](#) is among the 3% of the best Universities in the world, the Department of Mathematics is quite new (founded in 2019). First-year students have only attended online courses, while second-year have attended one full semester face-to-face (only the winter semester 2019-2020) and three semesters online. Another important element is that the Department located in central Greece, is easily accessible and chosen by young people of any cultural, economic and social background from all over the country. Students have the opportunity to take many courses, which combine theoretical foundation with applications in the context of the *interdisciplinary approach* (Rizos 2018). Some of the courses that students attend, and in which we are interested for the purposes of our research, are the following:

- «Calculus» (real numbers, sequences, series, functions, limits, derivatives and integrals with many applications in Science and Technology)
- «Philosophy of Science» (philosophical currents such as empiricism, rationalism and induction, and views of philosophers like Kuhn, Lakatos, Popper and Feyerabend)
- «Physics» (classical Mechanics, vibrations, waves and thermodynamics)
- «Programming Languages» (Python, MATLAB etc., laboratory exercises) and
- «Digital Technologies in Mathematics Education» (learning theories, integration of technology and science instruction into everyday math classroom, digital scenarios, dynamic geometry environments, didactical and social consequences).

Taking into account the current educational and social reality, we were interested to find out the ways students reacted regarding learning Sciences in a distance learning environment. Below we present an e-questionnaire we designed and we discuss the answers given by the students. The aim is to elaborate on the level and the ease of adaptability of the students to the remote processes. Our initial assessment was that the distribution of people's thesis on adaptability to change can be shown in a Gauss curve (Rogers 1962, p. 247). On the far left side of the curve there are the innovators and on the far right side there are the laggards. Between those ends there are some more categories, like early adopters, early majority and late majority. Adaptability in education (Green et al. 2020) on individual and community levels is very significant, especially if one considers education as a *public good* as we do. It is therefore essential to know in what "section" of the above curve our students are and to listen to their suggestions attentively, in order to improve the teaching procedure and academic life.

## 4. Overview of the project

At the end of the spring semester 2020-2021 we conducted a research project on Greek math students' experiences of distance learning using ICT. We asked, via email, a number of students from our institution to fill in an anonymous e-questionnaire which was designed on the MS-Forms platform. The criterion for the selection of the students was to have taken Calculus, Physics, Philosophy of Science, Programming Languages and Digital Technologies in Mathematics Education courses in the 2020-2021 academic year. We collected answers from 50 different students, specifically 25 from each academic year.

The students were asked 14 multiple choice questions and 2 open-ended questions. The questions were related to the obstacles (technical, academic, personal) faced by students due to distance learning, the skills they developed to deal with some of the above challenges, the impact of the pandemic on their lives and their suggestions for the next day. Furthermore, we talked via MS-Teams with 20 out of the 50 students and asked for more information about the topics. In the discussion we had, they reported the use of social media and internet platforms in education of mathematics and in communication between them. Some of the apps used for communication were Messenger, Facebook, Viber, Instagram etc. In general, chat spaces and video call apps like Skype and Zoom were used throughout the year for communication with the lecturers, too. It is quite obvious that this type of apps were very useful to students during the pandemic. By analyzing their answers we found some very interesting data about how our students try to learn Mathematics and what really helps them understand the content of each subject and each lecture. It was mentioned that students watch multiple videos related to math and physics in order to understand the concept better and they constantly use modern digital technologies that help the learning procedure. They said that the app GeoGebra really helped them understand graphs and two dimensional geometrical shapes and polygons. Other apps, like WolframMath help them calculate series, integrals, derivatives and more. They told us that they use programming languages and platforms for further math applications, some of them were Octave, MATLAB, Maxima, LateX etc.

Due to the pandemic the universities in Greece tend to come through a digital transformation (cf. Iivary et al. 2020). This educational (and economic) *digital transformation* is affecting and will continue to affect many different aspects of our lives and it seems to be a sign of an upcoming *social transformation*. This can be seen in the analysis of the answers to the questionnaire that follows. From the above analysis, we can study all the ways of transformation in science, society and academic life.

**5. Analysis**

We got 50 answers in our questionnaire; 25 were on the second semester and the remaining 25 were on the fourth semester. At first, we asked the students how much they think they have adapted to the distance learning process. Over 76% answered that they have adapted "enough" or "a lot" (see Table I). It seems that most of our students have easily and quickly adapted to the new model of education. Then, we asked about the technical problems the students encountered during the past academic year. In their responses, it was evident that by far the most common technical issue was the poor Wi-Fi connection and multiple network errors. This issue appeared in more than 30 questionnaires – the students had the option of answering this question by reporting multiple issues. Regarding the parameters that caused problems to the students in the learning process (solving exercises, applications and projects) the answers were pretty much balanced out. Every choice that we set beforehand has been chosen more than 15 times. The parameters were the following: "no access to libraries", "no access to laboratories", "communication with teachers", "communication with students", "lack of references". After addressing the basic problems, we need to concentrate on what type of relationship is dominant in our academic environment.

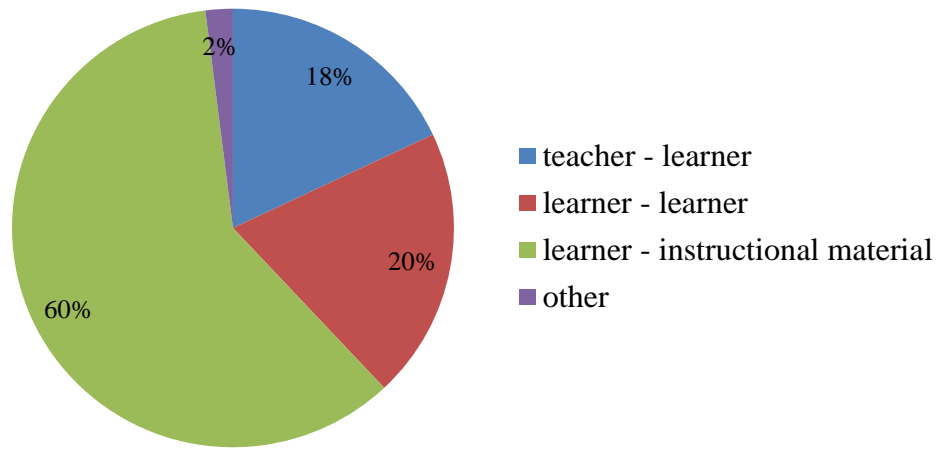

**The dominant relation in e-learning**

- teacher - learner (18%)
- learner - learner (20%)
- learner - instructional material (60%)
- other (2%)

The vast majority of the participants, over 60%, agree that the relationship between instructional material – learner is the projecting one (cf. Anderson 2008). During the e-learning era it seems that the students set the basis of their knowledge on books, papers and videos. The relationship with the instructional material seems to have an effect on the work load of the undergraduate students.

One of the most important things in e-learning and the digital transformation of the university is the use of modern platforms and apps. So, we asked our students about their developing knowledge in the use of those platforms. The majority of 62% reported that their skills have been upgraded "enough" or "a lot" (see Table I). The formalistic part of the experience in every university is the semester exams. Therefore, we wanted to know how our students think the exams would be more effective. We asked them in what model they would like the exams to be conducted. It seemed that the students wavered among different opinions. 40% of the participants preferred exams via the Internet, 38% of them preferred exams at the auditorium and the remaining 22% did not know what they preferred.

Table I

|  | Not at all | A little | Enough | A lot | Totally |
|---|---|---|---|---|---|
| How much do you believe you have adapted to the situation of distance learning in the University? | 0 (0%) | 6 (12%) | 23 (46%) | 15 (30%) | 6 (12%) |
| How much do you consider that your skills in using programs and apps (Python, GeoGebra, MATLAB etc.) during the academic year have been improved? | 6 (12%) | 13 (26%) | 19 (38%) | 11 (22%) | 1 (2%) |
| How much have your social life and your family equilibrium changed because of the pandemic? | 5 (10%) | 7 (14%) | 21 (42%) | 8 (16%) | 9 (18%) |

The pandemic has affected not only the part of the university but the aspect of the students' everyday lives, too. It is obvious that the modern economy has been severely damaged by the Covid-19 crisis, but how much have the economic changes affected Greek undergraduate students? When we asked them, 72% answered that the economic reverberations in their studies were "serious" or "computable". One fifth of the participants answered that their economic situation did not change at all. From their answers we can estimate the number of students who are in a difficult economic situation and the number of them who live in comfort. But, the pandemic has influenced not only the economic factor, but also the social lives and the family equilibrium of the students. After they were asked, 76% answered that their social lives have changed "enough", "a lot" or "totally" (see Table I).

Time management is a very important part in academic life and in a university in general. So, we asked in what type of lessons the students would spend more time on. In the multiple choice question, 44% of them answered that they would spend more time on face-to-face lessons, 14% that would spend more time on e-lessons, 34% answered would spend the same time on both and the remaining 8% did not know.

We specifically asked about the factors that caused difficulties to the students, offering them the option of multiple answers to the question. The most popular answer, by 46 votes, was the fatigue caused by the numerous hours in front of a wide screen. When the students were asked about the work load they had had during distance learning in comparison with face-to-face learning, 42% said that they would have lighter work load if the lessons were made face-to-face. The remaining percentage is split in other options: "the same load", "heavier load", "I do not know". This probably happens because in the online lessons students are called to solve more exercises for each subject and take part in more projects. Students say it is a bit tiring, but it is a more effective way of learning and being an active member in the scientific community.

At the end of the questionnaire, we decided to focus on the future. One of the final questions we asked the students was how many days they would prefer to keep having e-lessons after the end of the pandemic. The most popular answer was 5 days a week with 12 answers. This is a very interesting answer because through the questionnaire it is visible that students probably would like to return to the lessons in the auditorium and when they were asked the above question, things seemed quite confusing. It is possible that their adaptability to change has been compromised due to multiple challenges that they have faced during the past year. An interesting concept for future research would be about the change of adaptability of an undergraduate student (and their "movement" on the Gauss curve of adaptability). Although 24% of the students answered 5 days, 0, 1 and 2 days got 20% each.

Bar chart I

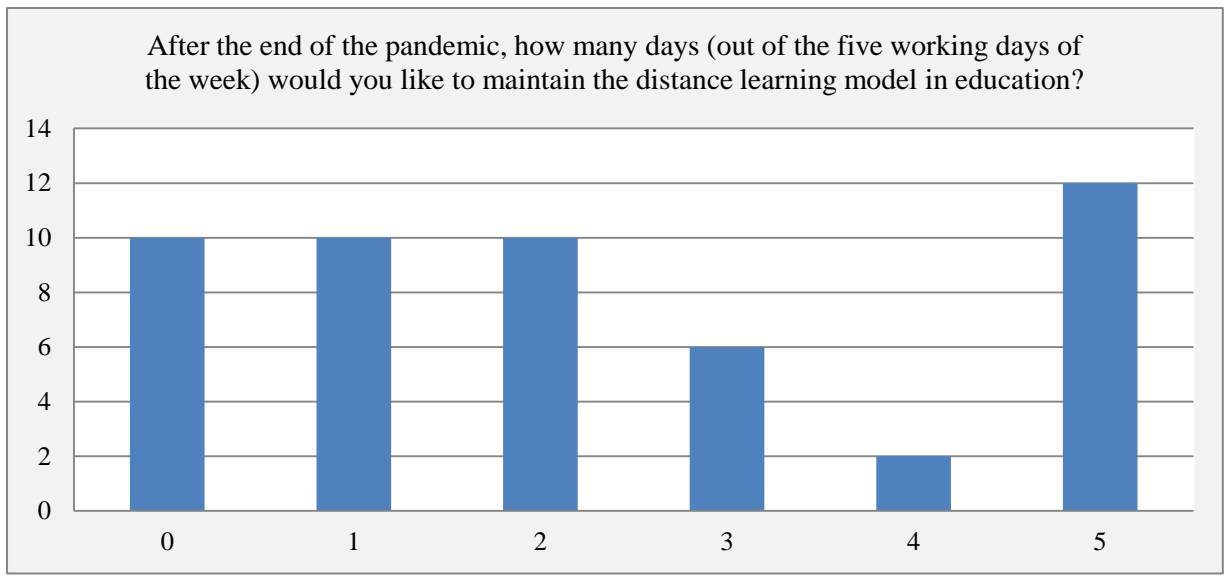

As after some years the students will be fresh professionals, the final multiple choice question was about the job market and, in similar form with the previous one, how many days of the week they would like to work from home, in their future career. The most popular answer was 2 days per week with 32% of the votes. These results are important for the digital transformation of the universities and the job market in the future.

Bar chart II

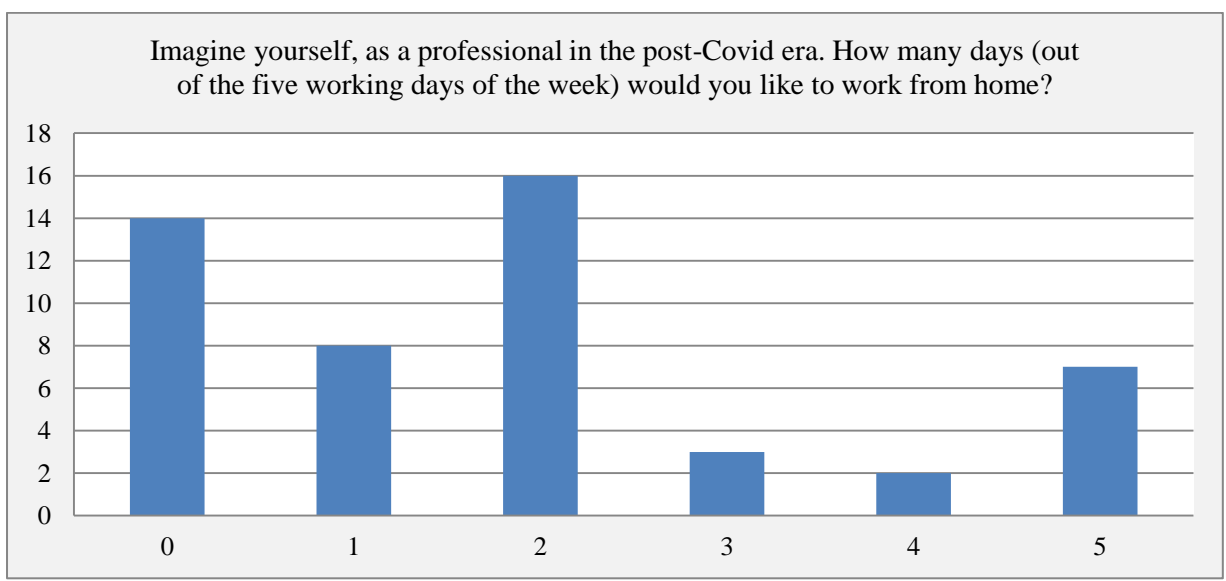

In addition, we also have included in the questionnaire two questions where students could answer freely in text. The first question was about their complete experience in the university during the academic year 2020-2021. In order to present the full experience of our students during the pandemic and their suggestions, we are going to quote some of the answers we got and some more interviews. We chose the most typical answers on each topic.

**Question 1:** *Describe your whole experience from distance learning in the Department of Mathematics of University of Thessaly during the academic year 2020-2021.*

**Student 1:** «*During this academic year, the learning process of Mathematics was quite difficult, because of the lack of interactivity. Another major issue was the time schedule and management. When nobody is used to e-learning, teachers cannot really give enough time to students in order to think about a question, which has an effect on the semester exams, too. The fact is that each side, professors and students, suit from distance learning. The teachers can complete a larger amount of syllabus without actually being understood by skipping through multiple slides and the students can get easily distracted and concentrate on other things without participating in the session. It seems that everyone has been tired of this situation and its functionality in Greek universities is questionable*».

**Student 2:** «*[…] but every teacher and student can conduct and attend a lesson at the comfort of their homes. People from different cities are not obliged to rent a house near the university or spend long hours using public transportation*».

**Student 3:** «*We have actually adapted to the new situation, apart from the social problems and lockdowns. It also depended on my psychological state during the past semesters*».

**Student 4:** «[…] *we learned to use the modern web platforms and we upgraded our computer skills, both personally and as a part of the academic community. But, we faced a problem in conducting team projects because of the lack of collaboration tools and communication. The communication problem was present in both teacher-student collaboration and student-student collaboration. The Internet connection did not help us either, because we faced multiple problems throughout each session*».

**Question 2:** *What are your suggestions for the next day of the pandemic in Education and the job market?*

**Student 5:** «*The truth is that we have been living in uncertainty for more than a full year and probably this uncertainty will continue to bedevil us even after the end of the pandemic. Of course, even after the pandemic we need to ensure a feeling of safety in the society and in education specifically. I think that the new semester should be executed with face-to-face learning because of the knowledge gaps created by poor performance and lack of interaction in distance learning*».

**Student 6:** «*One of my suggestions for the future in education is the use of a blended model of distance and face-to-face learning. A good combination for the students would be, to be allowed to attend theoretical lessons via Ms-Teams and applied lessons with physical presence. Another, probably functional way of using a blended model is assigning a specific number of days when the sessions would be conducted electronically. For example, two days of the week the lessons will be conducted with distance learning using ICT and the remaining three would be conducted face-to-face*».

Some of the problems the students encountered, as arising from their responses, are already known to all Greek universities and their presence is obvious in the academic community. By taking seriously into consideration students' responses and reshaping our educational model, we will try to overcome those problems.

On the suggestions topic, we read about a blended model, a flipped classroom model and a lot of other interesting ideas regarding education and economy. By integrating those ideas we aim, in addition to studying the upcoming digital and social transformation, at reflection and discussion on educational policy issues.

**6. Discussion**

When we analysed the responses to our questionnaire we found an interesting series of *paradoxes* when coupling pairs of answers. In fact, the validity of the existence of those paradoxes is accurate because these questions were asked in succession. The most significant and interesting paradox is about the already existing education model and the one students prefer for the next day. When students were asked if they liked face-to-face lessons or internet lessons, the vast majority answered that they would like to return to face-to-face lessons and in the next question they answered that they would like to keep the internet lessons for five days per week. We are almost confident that a factor that affects those answers is the adaptability of the students. Even if they can understand the reduction of the quality of the e-lessons, they would like to continue with it only because they do not want to change their everyday routine, their everlasting pattern.

Another paradox, noted in the answers we received, is the fact that our students like using digital technologies despite network errors. There seems to be an unreasonable "technophobia" in the student community, which might be because of the extensive use of digital technologies in their everyday lives, due to the lockdowns and even the rejection of modernity from their side. It is probable that a student sceptical to technology is experiencing his own "Groundhog day", which means that if the use of digital technologies was not extensive they would appreciate its abilities, but after using it approximately 8 hours a day it seems like a torture to them. This illusion is quite common in literature and everyday life, like the "torment of Sisyphus", a tragic figure of Greek mythology trapped in an everlasting loop of getting a boulder up a mountain and then rolling back down. So, this is a possible factor for the creation of the above paradox. Those "social" paradoxes are worth being studied in a future research.

Our students study Science and Technology in most of their sessions. Digital technologies are very useful in our profession and it is obvious that those illusions and paradoxes have to disappear in due time in order for them to concentrate on being productive and effective. In the next section, we are going to observe some of the uses of ICT in our Department and on STEM in general.

Students can use digital technologies in order to visualize theoretical problems and scenarios and to explore some real life applications (e.g. graphs, regular polygons, periodical movement simulations, representations of functions, fractals etc.). By using computers we can enhance the way Mathematics is taught and take the extra step. In this way both Pure and Applied Mathematics could evolve. By using web platforms and apps in Mathematics, students could probably understand better an exercise and maybe create an original teaching scenario themselves. «Digital Technologies in Mathematics Education» is a course that can possibly play the role of gateway for digital learning in mathematics education (Mulenga & Marbán 2020), widen the horizons of students and have a positive impact on their future careers (Borba 2021). From our conversation and interviews, we recognised that our students have multiple ideas for future research in Mathematics. On top of the above, we received the following ideas:

*«By using Python and the "numpy" and "matplot" libraries we can study multiple functions and even change graphic environment when we study a different geometry system»*. This student wants to specialize in Mathematics Programming.

*«By using GeoGebra we can create teaching scenarios and digital teaching scenarios in order to enhance the learning experience for every student. By using stories and graphic effects we can increase the students' engagement and help them understand theoretical notions and concepts»*. This student wants to specialize in Theoretical Mathematics and Mathematics Education.

Another topic that we talked about in our Department is the nature of this crisis. Mathematics, as a Science can format a crisis and help society overcome it (Skovsmose 2021). The most obvious crisis in this era is the Covid-19 outbreak and its impact in education. In 2019, Ole Skovsmose published a paper about the foundation of the relationships between Mathematics and crises, in which he concentrated on the flawed systems and the mathematical model behind educational system and the possibility of improving it (Skovsmose 2019). Our students' suggestions aim at that perspective of improvement in the system and mathematically reshaping and reconstructing a more efficient and complete structure for the institution.

We hope that this research will have an impact on our Department and Education itself. We attempted to present the reality in the academic society and to state that we should be open to new ideas and unconventional suggestions from young people. That way, we can improve ourselves and create a more stable, friendly and effective educational system. The future is ahead and we have the chance to shape it as we like. The point is that even when a crisis appears, we should work to get the best out of it and evolve.

**To cite this article:**